\documentclass[11pt]{article}

\usepackage{amssymb}

\newcommand{\newc}{\newcommand}

\newc{\eqnoset}{\setcounter{equation}{0}}

\newcommand{\mref}[1]{(\ref{#1})}
\newcommand{\reflemm}[1]{Lemma~\ref{#1}}
\newcommand{\refrem}[1]{Remark~\ref{#1}}
\newcommand{\reftheo}[1]{Theorem~\ref{#1}}

\newcommand{\refcoro}[1]{Corollary~\ref{#1}}

\newcommand{\refsec}[1]{Section~\ref{#1}}

\newcommand{\beq}{\begin{equation}}
\newcommand{\eeq}{\end{equation}}
\newcommand{\beqno}[1]{\begin{equation}\label{#1}}

\newcommand{\barr}{\begin{array}}
\newcommand{\earr}{\end{array}}

\newc{\bearr}{\begin{eqnarray*}}
\newc{\eearr}{\end{eqnarray*}}

\newc{\bearrno}[1]{\begin{eqnarray}\label{#1}}
\newc{\eearrno}{\end{eqnarray}}

\newc{\non}{\nonumber}
\newc{\nol}{\nonumber\nl}

\newcommand{\bdes}{\begin{description}}
\newcommand{\edes}{\end{description}}
\newc{\benu}{\begin{enumerate}}
\newc{\eenu}{\end{enumerate}}
\newc{\btab}{\begin{tabular}}
\newc{\etab}{\end{tabular}}

\newtheorem{theorem}{Theorem}[section]
\newtheorem{defi}[theorem]{Definition}
\newtheorem{lemma}[theorem]{Lemma}
\newtheorem{rem}[theorem]{Remark}
\newtheorem{exam}[theorem]{Example}
\newtheorem{propo}[theorem]{Proposition}
\newtheorem{corol}[theorem]{Corollary}

\newcommand{\btheo}[1]{\begin{theorem}\label{#1}}
\newc{\brem}[1]{\begin{rem}\label{#1}\em}
\newc{\bexam}[1]{\begin{exam}\label{#1}\em}
\newc{\bdefi}[1]{\begin{defi}\label{#1}}
\newcommand{\blemm}[1]{\begin{lemma}\label{#1}}
\newcommand{\bprop}[1]{\begin{propo}\label{#1}}
\newcommand{\bcoro}[1]{\begin{corol}\label{#1}}
\newcommand{\etheo}{\end{theorem}}
\newcommand{\elemm}{\end{lemma}}
\newcommand{\eprop}{\end{propo}}
\newcommand{\ecoro}{\end{corol}}
\newc{\erem}{\end{rem}}
\newc{\eexam}{\end{exam}}
\newc{\edefi}{\end{defi}}

\newc{\rmk}[1]{{\bf REMARK #1: }}
\newc{\DN}[1]{{\bf DEFINITION #1: }}

\newcommand{\bproof}{{\bf Proof:~~}}
\newc{\eproof}{{\vrule height8pt width5pt depth0pt}\vspace{3mm}}

\newc{\bfrac}[2]{\dspl{\frac{#1}{#2}}}

\newc{\nid}{\noindent}

\newcommand{\dspl}{\displaystyle}
\newc{\grad}{\nabla}
\newc{\Div}{\mbox{div}}
\newc{\pdt}[1]{\dspl{\frac{\partial{#1}}{\partial t}}}
\newc{\pdn}[1]{\dspl{\frac{\partial{#1}}{\partial \nu}}}
\newc{\pdNi}[1]{\dspl{\frac{\partial{#1}}{\partial \mathcal{N}_i}}}
\newc{\pD}[2]{\dspl{\frac{\partial{#1}}{\partial #2}}}
\newc{\dt}{\dspl{\frac{d}{dt}}}
\newc{\bdry}[1]{\mbox{$\partial #1$}}
\newc{\sgn}{\mbox{sign}}

\newc{\Hess}[1]{\frac{\partial^2 #1}{\pdh z_i \pdh z_j}}
\newc{\hess}[1]{\partial^2 #1/\pdh z_i \pdh z_j}

\newc{\ag}{\alpha}
\newc{\bg}{\beta}
\newc{\cg}{\gamma}\newc{\Cg}{\Gamma}
\newc{\dg}{\delta}\newc{\Dg}{\Delta}
\newc{\eg}{\varepsilon}
\newc{\zg}{\zeta}
\newc{\thg}{\theta}
\newc{\llg}{\lambda}\newc{\LLg}{\Lambda}
\newc{\kg}{\kappa}
\newc{\rg}{\rho}
\newc{\sg}{\sigma}\newc{\Sg}{\Sigma}
\newc{\tg}{\tau}
\newc{\fg}{\phi}\newc{\Fg}{\Phi}
\newc{\vfg}{\varphi}
\newc{\og}{\omega}\newc{\Og}{\Omega}
\newc{\pdh}{\partial}

\newc{\ccG}{{\cal G}}

\newc{\ii}[1]{\int_{#1}}
\newc{\iidx}[2]{{\dspl\int_{#1}~#2~dx}}
\newc{\bii}[1]{{\dspl \ii{#1} }}
\newc{\biii}[2]{{\dspl \iii{#1}{#2} }}
\newc{\su}[2]{\sum_{#1}^{#2}}
\newc{\bsu}[2]{{\dspl \su{#1}{#2} }}

\newc{\biiom}[1]{{\dspl\int_{\bdrom}~ #1 ~d\sg}}
\newc{\io}[1]{{\dspl\int_{\Og}~ #1 ~dx}}
\newc{\bio}[1]{{\dspl\int_{\bdrom}~ #1 ~d\sg}}
\newc{\bsir}{\bsu{i=1}{r}}
\newc{\bsim}{\bsu{i=1}{m}}

\newc{\iibr}[2]{\iidx{\bprw{#1}}{#2}}
\newc{\Intbr}[1]{\iibr{R}{#1}}
\newc{\intbr}[1]{\iibr{\rg}{#1}}
\newc{\intt}[3]{\int_{#1}^{#2}\int_\Og~#3~dxdt}

\newc{\itQ}[2]{\dspl{\int\hspace{-2.5mm}\int_{#1}~#2~dz}}
\newc{\mitQ}[2]{\dspl{\rule[1mm]{4mm}{.3mm}\hspace{-5.3mm}\int\hspace{-2.5mm}\int_{#1}~#2~dz}}
\newc{\mitQQ}[3]{\dspl{\rule[1mm]{4mm}{.3mm}\hspace{-5.3mm}\int\hspace{-2.5mm}\int_{#1}~#2~#3}}

\newc{\mitx}[2]{\dspl{\rule[1mm]{3mm}{.3mm}\hspace{-4mm}\int_{#1}~#2~dx}}
\newc{\mitmu}[2]{\dspl{\rule[1mm]{3mm}{.3mm}\hspace{-4mm}\int_{#1}~#2~d\mu}}
\newc{\iidmu}[2]{{\dspl\int_{#1}~#2~d\mu}}

\newc{\iidm}[3]{{\dspl\int_{#1}~#2~d #3}}

\newc{\itQmu}[2]{\dspl{\int\hspace{-2.5mm}\int_{#1}~#2~d\mu}}
\newc{\mitQmu}[2]{\dspl{\rule[1mm]{4mm}{.3mm}\hspace{-5.3mm}\int\hspace{-2.5mm}\int_{#1}~#2~d\mu}}

\newc{\mitQq}[2]{\dspl{\rule[1mm]{4mm}{.3mm}\hspace{-5.3mm}\int\hspace{-2.5mm}\int_{#1}~#2~d\bar{z}}}
\newc{\itQq}[2]{\dspl{\int\hspace{-2.5mm}\int_{#1}~#2~d\bar{z}}}

\newc{\pder}[2]{\dspl{\frac{\partial #1}{\partial #2}}}

\newc{\bdrom}{\bdry{\Og}}

\newc{\bilhom}{\mbox{Bil}(\mbox{Hom}(\RR^{nm},\RR^{nm}))}
\newc{\VV}[1]{{V(Q_{#1})}}

\newc{\ccA}{{\mathcal A}}
\newc{\ccB}{{\mathcal B}}
\newc{\ccC}{{\mathcal C}}
\newc{\ccD}{{\mathcal D}}
\newc{\ccE}{{\mathcal E}}
\newc{\ccH}{\mathcal{H}}
\newc{\ccF}{\mathcal{F}}
\newc{\ccI}{{\mathcal I}}
\newc{\ccJ}{{\mathcal J}}
\newc{\ccK}{{\mathcal K}}
\newc{\ccP}{{\mathcal P}}
\newc{\ccQ}{{\mathcal Q}}
\newc{\ccR}{{\mathcal R}}
\newc{\ccS}{{\mathcal S}}
\newc{\ccT}{{\mathcal T}}
\newc{\ccX}{{\mathcal X}}
\newc{\ccY}{{\mathcal Y}}
\newc{\ccZ}{{\mathcal Z}}

\newc{\bb}[1]{{\mathbf #1}}

\newc{\myprod}[1]{\langle #1 \rangle}
\newc{\mypar}[1]{\left( #1 \right)}

\newc{\BLLg}{\mathbf{\LLg}}

\newc{\mA}{\mathbf{A}}
\newc{\mB}{\mathbf{B}}
\newc{\mC}{\mathbf{C}}
\newc{\mD}{\mathbf{D}}
\newc{\mE}{\mathbf{E}}
\newc{\mF}{\mathbf{F}}
\newc{\mJ}{\mathbf{J}}
\newc{\mG}{\mathbf{G}}
\newc{\mP}{\mathbf{P}}
\newc{\mR}{\mathbf{R}}
\newc{\mQ}{\mathbf{Q}}
\newc{\mX}{\mathbf{X}}
\newc{\muu}{\mathbf{u}}
\newc{\mvv}{\mathbf{v}}

\newc{\mllg}{\mathbb{\lambda}}
\newc{\mLLg}{\mathbf{\LLg}}

\newc{\lspn}[2]{\mbox{$\| #1\|_{\Lsp{#2}}$}}
\newc{\Lpn}[2]{\mbox{$\| #1\|_{#2}$}}
\newc{\Hn}[1]{\mbox{$\| #1\|_{H^1(\Og)}$}}

\newc{\mynorm}[2]{\| #1\|_{#2}}

\newcommand{\RR}{{\rm I\kern -1.6pt{\rm R}}}

\newc{\itQQ}[2]{\dspl{\int_{#1}#2\,dz}}
\newc{\mmitQQ}[2]{\dspl{\rule[1mm]{4mm}{.3mm}\hspace{-4.3mm}\int_{#1}~#2~dz}}
\newc{\MmitQQ}[2]{\dspl{\rule[1mm]{4mm}{.3mm}\hspace{-4.3mm}\int_{#1}~#2~d\mu}}

\newc{\MUmitQQ}[3]{\dspl{\rule[1mm]{4mm}{.3mm}\hspace{-4.3mm}\int_{#1}~#2~d#3}}
\newc{\MUitQQ}[3]{\dspl{\int_{#1}~#2~d#3}}

\newc{\mccP}{\mathbb{P}}
\newc{\mccK}{\mathbb{K}}

\newc{\DKTmU}{\mccK(U)}
\newc{\DKTmUold}{(K_U(U)^{-1})^T}

\newc{\myPi}{\mathbf{W}}
\newc{\myIbar}{\bar{\ccI}_1}
\newc{\myIhat}{\hat{\ccI}_1}
\newc{\myIbreve}{\breve{\ccI}_0}

\newc{\mmk}{\mathbf{k}}

\newcommand{\ma}{\mathbf{a}}
\newcommand{\mg}{\mathbf{g}}

\newc{\mfu}{\mathbf{f_u}}
\newc{\mh}{\mathbf{h}}
\newc{\mb}{\mathbf{b}}

\newc{\mN}{\mathbf{N}}
\newc{\mI}{\mathbf{I}}
\newc{\mH}{\mathbf{H}}

\newc{\mk}{\mathbf{k}}
\newc{\mr}{\mathbf{r}}

\newc{\DIAGM}[2]{\left[\barr{ccc}#1&0\ldots&0\\
	\vdots&\ddots&\vdots\\0&\ldots0&#2\earr \right]}
\newc{\DiagM}[2]{\mbox{diag}\left[#1
	\cdots #2 \right]}
\newc{\vVEC}[2]{\left[\barr{c}#1\\
	\vdots\\#2\earr \right]}
\newc{\hVEC}[2]{\left[#1
	\cdots #2 \right]}

\newc{\mq}{\mathbf{q}}

\newc{\msys}[1]{\left\{\barr{l}#1\earr
	\right.}
\newc{\msysa}[1]{\left\{\barr{ll}#1\earr
	\right.}

\newc{\bbM}{\mathbb{M}}
\newc{\mat}[1]{\left[\barr{cc}#1\earr\right]}

\newc{\me}{\mathbf{e}}

\begin{document}

\vspace*{-.8in}
\begin{center} {\LARGE\em On the Continuity of Bounded Weak Solutions to  Parabolic Equations and Systems with Quadratic Growth in Gradients.}

 \end{center}

\vspace{.1in}

\begin{center}

{\sc Dung Le}{\footnote {Department of Mathematics, University of
Texas at San
Antonio, One UTSA Circle, San Antonio, TX 78249. {\tt Email: Dung.Le@utsa.edu}
}}

\end{center}

\begin{abstract}
We establish the pointwise continuity of bounded weak solutions to of a class of scalar parabolic equations and  strongly coupled parabolic  systems. Our approach to the regularity theory of parabolic scalar equations is quite elementary and its applications to strongly coupled systems does not require higher $L^p$ integrability of derivatives.   \end{abstract}

\vspace{.2in}

\section{Introduction} \label{intro}\eqnoset

Let $\Og$ be a bounded domain in $\RR^N$, $N\ge2$, with smooth boundary $\partial \Og$ and $T$ be a positive number. Denote $Q=\Og\times(0,T)$ and $z=(x,t)$ a generic point in $Q$. 

In the first part of this paper, we consider the following scalar parabolic equation
\beqno{0eq1}v_t=\Div(\ma Dv+\mb(v))+\me(Dv)+\mg(v) \quad \mbox{in $Q$}.\eeq
Here, $\ma,\me,\mg$ are scalar functions and $\mb\in\RR^N$. The equation is regular parabolic in the sense that $\ma$ is bounded and $\ma\ge\llg_0>0$ for some constant $\llg_0$.

Under suitable integrability conditions, much weaker than those in literature (e.g. \cite{LSU,Lieb}),  on the data of this equation we will show that any bounded weak solution $v$ of \mref{0eq1} is pointwise (or H\"older) continuous. Also important, we allow $\me(Dv)$ to have a quadratic growth in $Dv$. That is,   $\me(Dv)\le C(|Dv|^2+1)$.

In the second part, we will apply the theory for scalar equations to systems of $m$ equations, $m\ge2$, with linear or quadratic growth in gradients
\beqno{fullsys0aa}u_t=\Div(\mA(u)Du)+\me(Du)+f(u)
\eeq
in $Q$. Here, $u=[u_i]_{i=1}^m$ and $\ma(u)$ is a $m\times m$ matrix and $\me,f$ are vectors in $\RR^m$. We will always assume that there is $\llg_0>0$ such that for all $u\in\RR^m$, $\zeta\in\RR^{mN}$ and $i=1,\ldots,m$
\beqno{ellcondsys0} \myprod{\mA(u)\zeta,\zeta}\ge \llg_0|\zeta|^2.\eeq

Interestingly, we are able to establish the regularity of bounded weak solutions to parabolic systems on {\em planar domains}, where the 'hole filling' technique of Widman \cite{wid} for elliptic systems (e.g. \cite{BF}) does not seem to be extendable.
We also consider {\em triangular systems} on {\em any dimension} domains and assert that bounded weak solutions are pointwise continuity continuous. Our methods  does not require higher $L^p$ integrability of derivatives as in classical works (e.g. \cite{GiaS}).

In \refsec{techlem} we recall a simple parabolic version of the usual Sobolev inequality. We discuss scalar equations in \refsec{scalareqn}. We conclude the paper with applications to systems in \refsec{syseqn}.

\section{Some technical lemmas}\label{techlem}\eqnoset

In this section, we will present some technical results which will be used throughout this paper.
We recall the following simple parabolic version of the usual  Sobolev inequality 
\blemm{parasobo} Let $r=2/N$ if $N>2$ and $r\in(0,1)$ if $N\le 2$. If $g,G$ are sufficiently smooth then 
$$\itQ{\Og\times I}{|g|^{2r}|G|^2}\le C\sup_I\left(\iidx{\Og}{|g|^2}\right)^r\left(\itQ{\Og\times I}{(|DG|^2+|G|^2)}\right).
$$

If $G=0$ on $\partial\Og$ then we can drop the integrand $|G|^2$ on the right hand side. 

In particular, if $g=G$ we have
$$\itQ{\Og\times I}{|g|^{2(1+r)}}\le C(N,|\Og|)\sup_I\left(\iidx{\Og}{|g|^2}\right)^r\left(\itQ{\Og\times I}{(|Dg|^2+|g|^2)}\right).
$$

\elemm

Putting $g=G=|u|^p$ (with $p>1$) and using Young's inequality, we can see that for some  constant $c_0>0$ 
\beqno{pSobo}\left(\itQ{Q}{|u|^{2p(1+\frac 2N)}}\right)^\frac{N}{N+2}\le c_0\left(\sup_{(0,T)}\iidx{\Og}{|u|^{2p}}+\itQ{Q}{|u|^{2p-2}|Du|^2}\right).\eeq

In addition, if $\cg\in(1,1+2/N)$, then an use of Young and Sobolevs inequalities also gives that for any $\eg>0$ there is a constant $C(\eg)$ such that
\beqno{pSobo1}\|u\|_{L^{2p\cg}(Q)}\le \eg\left(\sup_{(0,T)}\iidx{\Og}{|u|^{2p}}+\itQ{Q}{|u|^{2p-2}|Du|^2}\right)^\frac{1}{2p}+C(\eg)\|u\|_{L^{2p}(Q)}.\eeq

\section{On scalar equations} \label{scalareqn}\eqnoset

\newc{\cM}{{\cal M}}

We now revisit the regularity theory of scalar equations with integrable coefficients, in the class $\cM(\Og,T)$ defined below. These new improvements  serve well our purposes in the next section.

For any $x_0\in\Og$, $R>0$ and $t_0\ge 4R^2$, we define $\Og_R(x_0)=\Og\cap B_R(x_0)$ and $Q_R(x_0)=\Og_R(x_0)\times(t_0-R^2,t_0)$. If $x_0,t_0$ are understood from the context, we simply drop them from the notations.

{\bf Definition of the class $\cM$:} We say that a function $f:Q\to \RR$ (or $\RR^m$) is of class $\cM(\Og,T)$ if for any $\eg>0$ there is $R(\eg)>0$ such that $\forall R\in(0,R(\eg))$ either
 $$\mathbf{i)}\quad \sup_{(0,T)}\|f\|_{L^{\frac{N}{2}}(\Og_R)} <\eg,$$
  or $$\mathbf{ii)}\quad \|f\|_{L^{\frac{N+2}{2}}(Q_R)} <\eg.$$

Alternatively, we also define the class $\bbM(\Og,T)$

{\bf Definition of the class $\bbM$:} We say that a function $f:Q\to \RR$ (or $\RR^m$) is of class $\bbM(\Og,T)$ if for  some $p_0>N/2$ such that  either one of the quantities i) $\sup_{(0,T)}\|f\|_{L^{p_0}(\Og)}$
or ii) $\|f\|_{L^{p_0+1}(Q)}$ is finite.

By H\"older inequality, it is is easy to see that $\bbM(\Og,T)\subset\cM(\Og,T)$.

\subsection{Global boundedness and a local estimate}
We consider scalar equation
\beqno{diagSKT} \msysa{v_t=\Div(ADv +B)+G& \mbox{in $Q$,}\\ v=v_0 & \mbox{in $\Og$.}} \eeq
Here  $A,G$ are scalar functions.
As usual, we will assume that there is a positive number $\llg_0$ such that
\beqno{ellcond} A\in\cM(\Og,T) \mbox{ and }A\ge \llg_0.\eeq
We also assume that there is a  function $\Fg\in\cM(\Og,T)$ such that $\Fg\ge\llg_0$ on $Q$ and
\beqno{grate}  |B|^2, |G| \le \Fg.\eeq

By using Steklov average, a weak solution of  \mref{diagSKT} satisfies for all $\eta\in C^1(Q)$
\beqno{weakeqn} \itQ{Q}{u_t\eta+ADvD\eta}=\itQ{Q}{G\eta}.\eeq

Applying the usual Moser iteration argument, we derive
 
\blemm{West} Assume \mref{ellcond} and the growth condition \mref{grate}  and $v$ is a solution of \mref{weakeqn}. Then  there is a constant $C$  such that \beqno{keyWesta}\sup_{Q}|v|\le C \left(\itQ{Q}{v^2}\right)^\frac12.\eeq

\elemm

\newc{\itQbar}[2]{\dspl{\int\hspace{-2.5mm}\int_{#1}~#2~d\bar{z}}}

The proof of this lemma bases on a Moser iteration technique by testing the equation with $|v|^{2p-2}v$ similar to the local version below.

We discuss the local estimates. This type of estimates will be useful for later investigations on the regularity of weak solutions. We will assume that the function $A$  is bounded. Note that $A$ may depend on $v$ in general. But $|v|$ is globally bounded by the above lemma.

\blemm{Westloc} Assume \mref{ellcond}, the growth condition \mref{grate} and that $A$  is bounded. Let $v$ be a solution of \mref{weakeqn} and $B_R$ be a ball in $\RR^N$. Then  there is a constant $C$ such that \beqno{keyWest}\sup_{Q_R}|v|\le C\left( \frac{1}{R^{N+2}}\itQ{Q_{2R}}{v^2}\right)^\frac12.\eeq

\elemm

We will make use of the Sobolev inequality for any $q\in(1,2N/(N-2))$ and $\eg>0$ there is $C(\eg)$ such that
\beqno{Sobo1}\left(\iidx{\Og}{|v|^q}\right)^\frac{2}{q}\le \eg\iidx{\Og}{|Dv|^2}+C(\eg)\iidx{\Og_R}{|v|^2}.\eeq

If $q=2N/(N-2)$ then we can only assert that 
\beqno{Sobo1a}\left(\iidx{\Og}{|v|^{\frac{2N}{N-2}}}\right)^\frac{N-2}{N}\le C\iidx{\Og}{|Dv|^2}+C\iidx{\Og_R}{|v|^2}.\eeq

The proof is the standard Moser iteration argument by testing \mref{weakeqn} with $|v|^{2p-2}v\fg^2\eta$ with some $p\ge1$ and $\fg,\eta$ are respectively  cutoff functions for concentric balls $B_R,B_{2R}$ and intervals $[-2R^2,-R^2]$, $[-R^2,0]$ with $|D\fg|\le C/R$, $|D\eta|\le C/R^2$. Let $V=|v|^{2p}$. Using \mref{grate} and Young's inequality  it is standard to derive (see \cite{Lieb,dleJMAA})
$$\sup_{(0,T)}\iidx{\Og}{|V|^2\fg^2\eta^2}+\itQ{Q}{|DV|^2\fg^2\eta^2}\le C\itQ{Q}{\Fg|V|^2\fg^2\eta^2}+\frac{C}{R^2}\itQ{Q}{|V|^2\fg^2\eta^2}.$$

By H\"older inequality with $q=2N/(N-2)$ so that $q/(q-2)=N/2$ and because $\Fg\in\cM(\Og,T)$ we see that, assuming i) in the definition of $\cM$,  for any $\eg>0$ if $R$ is sufficiently small then
$$ \barr{lll}\itQ{Q_{2R}}{\Fg|V|^2\fg^2\eta^2}&\le& \dspl{\int}_{-2R^2}^{0}{\left(\iidx{\Og_R}{\Fg^\frac{q}{q-2}}\right)^{1-\frac{2}{q}}\left(\iidx{\Og}{|V\fg\eta|^q}\right)^{\frac{2}{q}}}dt\\ &\le& \eg\dspl{\int}_{-2R^2}^{0}{\left(\iidx{\Og}{|V\fg\eta|^q}\right)^{\frac{2}{q}}}dt.\earr$$

Because $|D\fg|\le C/R$, we derive 
$$\eg\left(\iidx{\Og}{|V\fg\eta|^q}\right)^\frac{2}{q}\le \eg\iidx{\Og}{|DV|^2\fg^2\eta^2}+\frac{C(\eg)}{R^2}\iidx{\Og_{2R}}{|V|^2}.$$

If $\eg$ is sufficiently small in terms of $\llg_0$ then it follows that
\beqno{moser1}\sup_{(-R^2,0)}\iidx{\Og_R}{|V|^2}+\itQ{Q_R}{|DV|^2}\le \frac{C}{R^2}\itQ{Q_{2R}}{|V|^2}.\eeq

By the parabolic Sobolev inequality for $\cg=1+2/N>1$,  we obtain for any $p\ge1$
$$\left(\itQ{Q_R}{|v|^{2p\cg}}\right)^\frac{1}{\cg}\le \frac{C(\eg)}{R^2}\itQ{Q_{2R}}{|v|^{2p}}.$$

A standard Moser iteration argument (e.g. \cite{Lieb}) implies the local estimate of the lemma.

If ii) in the definition of $\cM$ holds then for $\cg=1+2/N$ we can also use H\"older inequality ($\cg'=(N+2)/2$) to have for $R$ small
\beqno{wholder}\itQ{Q_{2R}}{\Fg V^2}\le \left(\itQ{Q_{2R}}{\Fg^{\cg'}}\right)^\frac{1}{\cg'}\left(\itQ{Q_{2R}}{V^{2\cg}}\right)^\frac{1}{\cg}\le\eg\left(\itQ{Q_{2R}}{V^{2\cg}}\right)^\frac{1}{\cg}.\eeq

Therefore, we can use the above estimate and the parabolic Sobolev inequality to treat the last integral to obtain the local estimate \mref{moser1} again and the proof can go on.

If $\Fg\in\cM(\Og,T)$, depending whether i) or ii) in its definition holds, we define a positive function  $\nu$ in $R\ge0$ by \beqno{nudef}\nu(R):=\sup_{(-R^2,0)}\left(\iidx{B_R}{\Fg^\frac{N}{2}}\right)^\frac{2}{N} \mbox{ or }\left(\itQ{Q_R}{\Fg^\frac{N+2}{2}}\right)^\frac{2}{N+2}.\eeq

It is clear that  $\nu$ is increasing and continuous at $0$ ($\nu(0)=0$). Also, because either
$$\itQ{Q_R}{\Fg}\le \int_{-R^2}^0 \left(\iidx{B_R}{\Fg^\frac{N}{2}}\right)^\frac{2}{N}(R^{N})^\frac{N-2}{N}dt\le \nu(R)R^N$$
or $$\itQ{Q_R}{\Fg}\le \left(\itQ{Q_R}{\Fg^\frac{N+2}{2}}\right)^\frac{2}{N+2}(R^{N+2})^\frac{N}{N+2}=\nu(R)R^N,$$
we observe that  in both cases \beqno{Fglocbound} \itQ{Q_R}{\Fg}\le \nu(R)R^N.\eeq

\subsection{H\"older continuity} 

We will study the H\"older regularity in this subsection. Note that the bounds for the H\"older norm and exponents will depend only on the generic constants  the parameters in their definition $\cM(\Og,T)$. This fact will play a crucial role when we estimate the derivatives which appear in cross diffusion systems.

\btheo{blemm} Assume that $\ma\ge\llg_0$  and $\ma$ is bounded and $|\mb|^2,|\mg|\le\Fg\in \cM(\Og,T)$. Let $v$ be a bounded weak solution of 
\beqno{eq1}v_t=\Div(\ma Dv+\mb)+\me(Dv)+\mg.\eeq
If $\me(Dv)\le \eg_*|Dv|^2 +\Fg$ for some $\eg_*>0$ with $\eg_*\sup_Q|v|$ is sufficiently small compared to $\llg_0$ then $v$ is pointwise continuous. The continuity of $v$ depends on those of the integrals in \mref{nudef} on the measure of their domains.
\etheo

For simplicity we assume first that $\mb\equiv0$. The case $\mb\ne0$ is similar will be discussed in \refrem{brem} after this proof.

The idea based on that of \cite{dleNA}. We present the details and nontrivial modification. 

Fixing any $x_0\in\Og$, $t_0>0$ and $4R^2<t_0$, we denote $Q_{iR}=\Og_{iR}\times(t_0-iR^2,t_0)$ where $\Og_R=\Og\cap B_R(x_0)$.

Let $M_i=\sup_{Q_{iR}}v$, $m_i=\inf_{Q_{iR}}v$ and $\og_i=M_i-m_i$. For some number $\nu_*>0$ will be determined later and the function $\nu$ as in \mref{Fglocbound} we define $\dg(R)=\nu_*\sqrt{\nu(R)}$ and
$$ N_1(v)=2(M_4-v)+\dg(R),\quad N_2(v)=2(v-m_4)+\dg(R),$$
$$w_1(v)=\log\left(\frac{\og_4+\dg(R)}{N_1(v)}\right),\quad w_2(v)=\log\left(\frac{\og_4+\dg(R)}{N_2(v)}\right).$$

We will prove that either $w_1$ or $w_2$ is bounded from above on $Q_{2R}$ by a constant $C$ {\em independent of} $R$. This implies a decay estimate for some $\eg\in(0,1)$ and all $R>0$ \beqno{decay}\og_2\le \eg \og_4 + C\dg(R). \eeq It is standard to iterate \mref{decay} to obtain  (see \cite[Lemma 8.23]{GT}) 
$$\og(R)\le C\left[\left(\frac{R}{R_0}\right)^\ag\og(R_0)+\dg(R^\mu R_0^{1-\mu})\right]\quad \forall R\in(0,R_0)$$
for any $\mu\in(0,1)$ and some $R_0,\ag>0$.
This gives the continuity of $v$ as $\lim_{R\to0}\dg(R)=0$. 
To see \mref{decay}, if either $w_1$ or $w_2$ is bounded from above by $C>0$ in $Q_{2R}$ then either $$\og_4+\dg(R)\le 2C(\og_4+m_4-v)+C\dg(R)\mbox{ or }\og_4+\dg(R)\le 2C(\og_4+v-M_4)+C\dg(R).$$ 
Taking the supremum (respectively infimum) over $Q_{2R}$ and replacing $m_4$ by $m_2$ (respectively $M_4$ by $M_2$), we obtain
$\og_2\le \eg \og_4+(C-1)\dg(R)$ for $\eg=\frac{2C-1}{2C}<1$. This yields \mref{decay}.

Thus,  we just need to show that either $w_1$ or $w_2$ is bounded from above on $Q_{2R}$.
Before proving this, we note the following crucial property of the functions $w_1,w_2$. We will see that $w_1\le 0\Leftrightarrow w_2\ge0$ and vice versa. Indeed,
\beqno{alt}w_1\le0 \Leftrightarrow \og_4+\dg(R)\le 2(M_4-v)+\dg(R) \Leftrightarrow 2(v-m_4)+\dg(R)\le \og_4+\dg(R)\Leftrightarrow w_2\ge 0.
\eeq

\bproof For any $\eta\in C^1(Q)$ and $\eta\ge0$, observe that $Dw_1=\frac{2Dv}{N_1(v)}$,  $(w_1)_t=\frac{2v_t}{N_1(v)}$
and $Dw_2=-\frac{2Dv}{N_2(v)}$,  $(w_2)_t=-\frac{2v_t}{N_2(v)}$. So, by multiplying the equation of $v$ by $\eta/N_1(v)$ and $-\eta/N_2(v)$ and writing $w_i, N_i(v)$ respectively by $w, N(v)$ we obtain
\beqno{eq2a} \iidx{\Og}{\frac{\partial w}{\partial t}\eta}+ \iidx{\Og}{\myprod{\ma Dw,D\eta}}+\iidx{\Og}{\myprod{\ma Dv,\frac{\eta Dv}{N^2(v)}}}=2\iidx{\Og}{\frac{\me(Dv)\pm \mg}{N(v)}\eta}.\eeq

Because $\myprod{\ma Dv, Dv}\ge\llg_0|Dv|^2$ and the assumption $\me(Dv)\le \eg_*|Dv|^2+\Fg$ $$\left|\frac{\me(Dv)}{N(v)}\right| \le N(v)\frac{\eg_*|Dv|^2+\Fg}{N^2(v)}\le \eg_*(4\sup_Q|v|+1)\frac{|Dv|^2}{N^2(v)}+\frac{\Fg}{N(v)},$$ assuming $\nu(R_0)<1$,  we can absorb (and discard) the integral of $\frac{\me(Dv)\eta}{N(v)}$ into that of $\myprod{\ma Dv,\eta Dv/N^2(v)}$ if $\eg_*\sup_Q|v|$ is sufficiently small (compared to $\llg_0$)  to  get
\beqno{eq2}\iidx{\Og}{\frac{\partial w}{\partial t}\eta}+ \iidx{\Og}{\myprod{\ma Dw,D\eta}}\le 4\iidx{\Og}{\frac{|\Fg|}{N(v)}\eta}.\eeq

Testing \mref{eq2} with $(w^+)^{2p-1}\eta^2$  as in \reflemm{Westloc} with $G=\Fg/N(v)$. Because $\llg_0\le \ma$ and $\ma$ is bounded and since $N(v)\ge \nu_*\nu(R)$ and by the definition \mref{nudef} of $\nu$ then either $$\|G\|_{L^\frac{N}{2}(\Og_R)}\le \frac{1}{\nu_*\nu(R)}\|\Fg\|_{L^\frac{N}{2}(\Og_R)}=\frac{1}{\nu_*}$$ or $\|G\|_{L^\frac{N+2}{2}(Q_R)}\le \frac{1}{\nu_*\nu(R)}\|\Fg\|_{L^\frac{N+2}{2}(Q_R)}\le \frac{1}{\nu_*}$. Thus,
we see that $G\in \cM(\Og,T)$ if $\nu_*$ large. By \reflemm{Westloc}, which applies if $1/\nu_*$ is sufficiently small, we  find a  constant $C$ such that
\beqno{supw}\sup_{\Og_{2R}\times(t_0-2R^2,t_0)} w^+\le C\left(\frac{1}{R^{N+2}}\itQ{\Og_{4R}\times(t_0-4R^2,t_0)}{w^2}\right)^\frac12.\eeq

If we can show that for any $R>0$ there is a constant $C$ such that
\beqno{wbound} \frac{1}{R^{N+2}}\itQ{\Og_{4R}\times(t_0-4R^2,t_0)}{w^2}\le C\eeq then this implies $w^+$ is bounded on $Q_{2R}$. So, the decay estimate \mref{decay} holds. 

Let $\eta$ be a cut-off function for $B_{2R}, B_{4R}$. Replacing $\eta$ in \mref{eq2a} by $\eta^2$ (keeping the third term on the left hand side and using the assumption that $\eg_*\sup_Q|v|$ is small again), we get
\beqno{k1}\frac{d}{dt}\iidx{\Og}{w\eta^2}+\iidx{\Og}{\ma |Dw|^2\eta^2}\le \iidx{\Og}{\ma|Dw|\eta|D\eta|}+4\iidx{\Og}{\frac{|\Fg|}{N(v)}\eta^2}
\eeq

Applying Young's inequality we derive (as $|D\eta|\le C/R$ and $\ma$ is bounded)
\beqno{k2}\barr{lll}\frac{d}{dt}\iidx{\Og}{w\eta^2}+\iidx{\Og}{\ma |Dw|^2\eta^2}&\le& \frac{1}{R^2}\iidx{\Og}{\ma}+4\iidx{\Og}{\frac{|\Fg|}{N(v)}\eta^2}\\&\le& CR^{N-2}+4\iidx{\Og}{\frac{|\Fg|}{N(v)}\eta^2}.\earr
\eeq

Set $I_*=[t_0-4R^2, t_0-2R^2]$, $Q_*=B_{2R}\times I_*$ and $Q_v=\{(x,t)\in Q_*\,:\,\; w_1\le0\}$. It is easy to see that $w_2\le0$ on $Q_*\setminus Q_v$ (see \mref{alt}). Therefore one of $w_1^+,w_2^+$ must vanish on a subset $Q^0$ of $Q_*$ with $|Q^0|\ge \frac12|Q_*|$. We denote by $w$ such function. Let $Q_t^0$ be the slice $Q^0\cap (B_{2R}\times\{t\})$ then $Q^0=\cup_{t\in I_*} Q_t^0$.  For $t\in I_*$ let
$$\Og^0_t=\{x\,:\, w^+(x,t)=0\}, \quad m(t)=|Q^0_t|.$$
The fact that $|Q^0|\ge \frac12|Q_*|$ implies $\int_{I_*}m(t)dt\ge \frac12 R^{N+2}$.

We now set $$V(t)=\frac{\iidx{\Og}{w\eta^2}}{\iidx{\Og}{\eta^2}}.$$

By the weighted Poincar\'e' inequality (\cite[Lemma 3]{Moser}) 
$$\iidx{\Og}{(w-V)^2\eta^2}\le CR^2\iidx{\Og}{|Dw|^2\eta^2}.$$
Reducing the integral on the left to the set $Q^0_t$ where $w\le0$ (so that $V^2\le (w-V)^2$), we have $$V^2(t)m(t)\le CR^2\iidx{\Og}{|Dw|^2\eta^2}.$$
Since $N(v)\ge \nu_*\nu(R)$ on $Q_{4R}$, the above estimate and \mref{k2} implies that ($V'$ denotes the $t$ derivative)
\beqno{k3}R^N V'(t)+\frac{1}{R^2}V^2(t)m(t)\le CR^{N-2}+\frac{4}{\nu_*\nu(R)}\iidx{\Og}{|\Fg|\eta^2}.\eeq

Because $\nu$ is given by \mref{nudef}, we also get from \mref{Fglocbound}
\beqno{mgest}\|G\|_{L^1(Q_R)}\le \frac{1}{\nu_*\nu(R)}\itQ{Q_R}{\Fg}\le CR^N.\eeq

We show that $V(t_1)$ is bounded on for some $t_1\in I_*$. Indeed, suppose $V(t)\ge A>0$ in $I_*$. We have from \mref{k3} 
$$ R^{N+2}\frac{V'(t)}{V^2(t)}+m(t)\le \frac{C}{A^2}\left(R^N+\frac{R^2}{\nu_*\nu(R)}\iidx{\Og}{|\Fg|\eta^2}\right).$$

Because $\int_{I_*}m(t)dt\ge \frac12 R^{N+2}$ and $|I_*|\sim R^2$, we integrate this over $I_*$ and use \mref{mgest} to see that
$$R^{N+2}\le \int_{I_*}m(t)dt\le R^{N+2}\left(\frac{2}{A}+\frac{C}{\nu_*A^2}\right).$$
By choosing $A$ large we get a contradiction. So, we must have $V(t_1)\le A$ for some $t_1\in I_*$.

Integrating \mref{k2} over $[t_1,t_2]$ for any $t_2\in I_0=[t_0-2R^2, t_0]$, we have
$$V(t_2)\iidx{\Og\times\{t_2\}}{\eta^2}+\itQ{\Og\times I_0}{\ma |Dw|^2\eta^2}\le CR^N +V(t_1)\iidx{\Og\times\{t_1\}}{\eta^2}.$$
This implies that $V(t_2)\le C$ for all $t_2\in I_0$ and $\itQ{\Og\times I_0}{\ma |Dw|^2\eta^2}\le CR^N$. This implies that $V(t_2)\le C$ for all $t_2\in I_0$ and $\itQ{\Og\times I_0}{\ma |Dw|^2\eta^2}\le CR^N$. Since we can always assume that $\og_4\ge \nu(R)$ (otherwise there is nothing to prove) so that $V(t)$ is bounded from below by $-\log(\nu_*)$. Hence, $|V(t)|$ is bounded.

By Poincar\'e's inequality again we have
$$\itQ{\Og\times I_0}{(w-V)^2\eta^2}\le CR^2\itQ{\Og\times I_0}{|Dw|^2\eta^2}\le CR^{N+2}.$$

Since $|V(t)|$ is bounded on $I_0$, replacing $R$ by $2R$, the above implies the desired \mref{wbound} and concludes our proof. \eproof

\brem{brem}  The assertion still holds if $\mb\ne0$ and $|\mb|^2\in\cM(\Og,T)$.
Indeed, there will be an extra term in our argument. Namely, replacing $\eta$ by $\eta/N(v)$ as before we have the following extra terms in \mref{eq2a}
$$\iidx{\Og}{\myprod{\mb, \frac{D\eta}{N(v)}+\frac{\eta Dv}{N^2(v)}}}=\iidx{\Og}{\myprod{\frac{\mb}{N(v)}, D\eta+\frac12 Dw \eta}}.$$

We take $\eta$ to be $|w|^{2p-2}w\eta^2$ and think of $B_1:= \mb/N(v), B_2:=\frac12\mb/N(v)$ in 
$$w_t=\Div(ADw+B_1)+B_2Dw+G$$
(compared with \mref{diagSKT}). Again noting that $Dw=2Dv/N(v)$. Because $|\mb|^2$, $|B_1|^2,|B_2|^2\le\Fg$  so that the Moser's technique is applied as before. We  get the local estimate \mref{supw} and the proof can go on (we can redefine $N(v)$ such that $N^2(v)\ge \nu_*\nu(R)$ because certainly  $\nu(R)^\frac{1}{2}\ge \nu(R)$).

\erem

\brem{FgVrem} It is worth to noting that if $\Fg$ verifies \mref{wholder} then \mref{supw} holds. Once this is true  we need only \mref{mgest} which is valid if the $L^1(Q)$ norm of $\Fg$ satisfies $\|\Fg\|_{L^1(Q_R)}\le \nu(R)R^N$ (see \mref{Fglocbound}). 

\erem

Since $\bbM(\Og,T)\subset \cM(\Og,T)$ a similar argument shows that we can take $\nu_*=1$ and $\nu(R)=R^\cg$ for some appropriate $\cg>0$ to have a stronger version of \reflemm{blemm}.

\blemm{bholderlemm} Assume that $\ma\ge\llg_0$ and $\ma$ is bounded and $|\mb|^2,|\mg|\le\Fg\in \bbM(\Og,T)$. Let $v$ be a bounded weak solution of 
$$v_t=\Div(\ma Dv+\mb)+\me(Dv)+\mg.$$
If $\me(Dv)\le \eg_*|Dv|^2 +\Fg$ for some $\eg_*>0\sup_Q|v|$ small compared to $\llg_0$ then $v$ is H\"older continuous. Its H\"older norm is bounded in terms of the $L^{p_0}(\Og)$ (or $L^{p_0+1}(Q)$) norms of $|\mb|^2,\mg$.
\elemm

It is also important to mention the following

\bcoro{Dvloc} Assume that $\ma=\ma(v)$ is H\"older continuous in $v$ and $\ma\ge\llg_0$ and that $|\mb|^2,\me$ are as in \reftheo{blemm}. Let $v$ be a bounded weak solution of 
$$v_t=\Div(\ma(v) Dv+\mb)+\me(Dv)+\mg.$$
If $\mg\in L^\infty_{loc}(Q)$ then $Dv$ is H\"older continuous and locally bounded.

\ecoro

\bproof Knowing that $v$  is continuous and so is $\ma(v)$, we can use \cite[Theorem 3.2]{GiaS} applying to scalar equations to see that $Dv$ is H\"older continuous and thus locally  bounded. Note that we don't have to use the $L^p$ estimate of $Dv$ here once we know that $v$ is continuous. This continuity of $v$ suffices to obtain \cite[(3.4) in the proof of Proposition 3.1]{GiaS} to obtain a decay estimate in proving that $Dv$ is H\"older continuous in the proof of \cite[Theorem 3.2]{GiaS}. \eproof

\section{Applications to systems} \label{syseqn}\eqnoset

In this section, we apply the theory to the system
\beqno{fullsys0a}u_t=\Div(\mA(u)Du)+\me(Du)+f(u)
\eeq
in $Q$. Here, $u=[u_i]_{i=1}^m$ and $\ma(u)$ is a $m\times m$ matrix and $\me,f$ are vectors in $\RR^m$. We will always assume that there is $\llg_0>0$ such that for all $u\in\RR^m$, $\zeta\in\RR^{mN}$ and $i=1,\ldots,m$
\beqno{ellcondsys} \myprod{\mA(u)\zeta,\zeta}\ge \llg_0|\zeta|^2  \mbox{ and }\ag_{ii}(u)\ge \llg_0.\eeq

In addition,  there is $\eg_*>0$ such that \beqno{egcond} |\me(\zeta)|\le \eg_*|\zeta|^2\quad \forall\zeta\in \RR^{mN}.\eeq

The well known 'hole-filling' trick of Widman (e.g. \cite{BF}) has been apllied in the regularity of strongly coupled {\em elliptic systems} on {\em planar} domains. Roughly speaking, the idea is that  one tests the elliptic system with $u\fg^2$ where $\fg$ is a cutoff function for $B_R,B_{2R}$ (that is $\fg\equiv 1$ in $B_R$ and $\fg\equiv0$ outside $B_{2R}$ and $|D\fg|\le C/R$) to obtain a decay estimate for $\|Du\|_{L^2(B_R)}$. This and imbedding theorems of Campanato spaces (e.g. \cite{Lieb}) implies that $u$ is H\"older continuous if $N=2$.

However, this trick does not seem to apply to the corresponding parabolic systems like \mref{fullsys0a} even when $\eg_*=0$ in \mref{egcond}. Following the same idea to \mref{fullsys0a} with $\fg$ is a cutoff function for $Q_R,Q_{2R}$, one can not obtain a decay estimate for $\|Du\|_{L^2(Q_R)}$. Even so, the extra time dimension does not imply any continuity of $u$.

We will apply the theory for scalar equations in previous section to each equation in \mref{fullsys0a} for appropriate conditions on  $\ma,\mb$ and $\mg$ to establish the pointwise continuity of bounded weak solutions to systems like \mref{fullsys0a}. Namely, under appropriate settings and assumptions, we will show that if a bounded weak solution $u$ of \mref{fullsys0a} is averagely continuous \beqno{avercont} \liminf_{R\to 0} \mitQ{Q_R}{|u-u_R|^2}=0\eeq
then it is pointwise continuous.

\subsection{Full systems (SKT) on planar domains} 
The checking of the average continuity assumption \mref{avercont} for a bounded weak solution to \mref{fullsys0a} is a hard problem in general. Here, we consider the case $N=2$ and present examples when this can be done.

Let $\mA=P_u$ for some $P:\RR^m\to\RR^m$. We consider a special case of \mref{fullsys0a}. The following model was introduced in \cite{SKT} and studied widely in the context of mathematical biology (e.g. \cite{yag}) \beqno{sktsys0a}u_t=\Delta(P(u))+\me(Du)+f(u)
\eeq
and assume that the nonlinearity is sublinear. That is, for some constant $C$ \mref{egcond} is now
\beqno{egcondsub} |\me(\zeta)|\le C|\zeta|\quad \forall\zeta\in \RR^{mN}.\eeq

We will prove that $\|Du\|_{L^2(\Og)}\le C$ for some constant $C$ for all $t\in(0,T)$. The following calculation is formal and it can be justified by replacing the operator $\frac{\partial}{\partial t}$ in the proof of \cite[Lemma 2.2]{dleJMAA} with the partial difference operator $\dg_h^{(t)}$ (or $u,P(u)$ by their Steklov average as in \cite{LSU}). Multiplying the $i^{th}$ equation with $\frac{\partial}{\partial t}P_i(u)\eta$ where $\eta$ is function in $t$ and summing the results, we obtain for $Q^t=\Og\times(t_1,t)$
$$\itQ{Q^t}{\myprod{u_t,P_u u_t}\eta}+\itQ{Q^t}{\myprod{D(P(u)), D(P(u))_t}\eta}=\itQ{Q^t}{\myprod{f,P_u u_t}\eta}.$$

We now choose $\eta$ such that $\eta(t)=1$, $\eta(t_1)=0$ and $|\eta'|\le C$. Since $\myprod{u_t,P_u u_t}\ge \llg_u|u_t|^2$, $|u|$ is bounded, by a simple use of Young's inequality to the last integral (assuming $f(u)\in L^2(Q)$ for any given bounded solution $u$) and rearranging, we easily derive ($C$ denotes a constant depending on $\sup_Q|u|$)
$$\iidx{\Og}{|Du(x,t)|^2}\le C\int_{t_1}^t\iidx{\Og}{|Du(x,t)|^2}dt +C.$$
This is an integral Gr\"onwall inequality for $y(t)=\|Du\|_{L^2(\Og\times\{t\})}$ and implies for all $t\in(0,T)$ that $\|Du\|_{L^2(\Og\times\{t\})}\le C$.

For each $i=1,\ldots,m$ we apply \reftheo{blemm} by simply set $\ma=P_{u_i}$, $\mb=\sum_{j\ne i}P_{u_j}Du_j$ and $\mg=\me(Du)+f$ (we see that $|\mb|^2,|\mg|$ belong to $\cM(\Og,T)$ as $N=2$ as either i) or ii) of the definition of $\cM$ is satisfied). Hence,  $u$ is pointwise continuous. 

\bcoro{SKT} Consider the system \mref{sktsys0a}. Assume that $N=2$, \mref{egcondsub} and $f\in L^2(Q)$. Then any bounded weak solution of \mref{sktsys0a} is pointwise continuous.

\ecoro

\brem{Dbounded} If $f\in L^\infty(Q)$ then the derivatives are bounded and H\"older continuous (see the discussion leading to \refcoro{trisysthm} below).
\erem

\subsection{Triangular systems on $N$-dimensional domains}

The result of \refcoro{SKT} holds for full systems with nonlinearities grow at most linear in gradients (see \mref{egcondsub}). If the system \mref{fullsys0a} is of the special triagular form then we can resume the quadratic growth in gradients \mref{egcond} (and some what more general) for general $N$.

We will present now an example of a class of  triangular systems whose nonlinearities having quadratic growth in gradients.

We start with a system of two equations
\beqno{trisys0}\msysa{u_t=\Div(\ag(u,v)Du+\bg(u,v)Dv)+\eg_1|Du|^2+c|Dv|^2+\mg_1&\mbox{in $Q$,}\\
	v_t=\Div(\dg(v)Dv)+\eg_2|Dv|^2+\mg_2&\mbox{in $Q$,}
	\\ \mbox{Homogeneous Dirichlet or Neumann boundary conditions}&\mbox{on $\partial\Og\times(0,T)$},
	\\
	(u,v)=(u_0,v_0)&\mbox{on $\Og$.}}\eeq

That is we will consider \mref{fullsys0a} with $\mA=\mat{\ag&\bg\\0&\dg}$, $\me(Du,Dv)=\left[\barr{c}\eg_{1}|Du|^2+c|Dv|^2\\\eg_{2}|Dv|^2\earr\right]$.

We assume, instead of \mref{ellcondsys} which implies,   that $\ag(u,v), \dg(v)\ge\llg_0$. We also assume that $\mg_i\in \cM(\Og,T)$ and $\ag,\bg$ are continuous and $\dg$ is H\"older continuous in $u,v$. 

Assume that $u,v$ are locally bounded. If $\eg_2\sup_Q|v|$ is small then we see that $v$ is continuous by \reftheo{blemm}, with $\mb=0$ and $\mg=\mg_2$.

Knowing that $\dg(v)$ is continuous and assuming that $\mg_1\in L^\infty_{loc}(Q)$, we can use \cite[Theorem 3.2]{GiaS} applying to scalar equations to see that $Dv$ is H\"older continuous and thus bounded (see \refcoro{Dvloc}). This can be used in the equation of $u$ with $\mb=\bg(u,v)Dv$ and $\mg=c|Dv|^2+\mg_2$ and we can apply \reftheo{blemm} again to prove that $u$ is continuous if $\eg_1\sup_Q|u|$ is small (and $c$ can be large). Again, note that we don't have to use the higher integrability $L^p$ estimate of $Dv$ here once we know that $v$ is continuous. This continuity of $v$ suffices to obtain \cite[(3.4) in the proof of Proposition 3.1]{GiaS} to obtain a decay estimate in proving that $Dv$ is H\"older continuous in the proof of \cite[Theorem 3.2]{GiaS}. 

Again, we remark that the coninuity of $u,v$ also shows that $Du,Dv$ are H\"older continuous by \cite{GiaS}. But to apply the theory in \cite{GiaS}, we need the elliptic condition \mref{ellcondsys} for the whole system when $m>2$.

By induction, the above argument can be extended to systems for $m>2$ unknowns $u_i$ ($i=1,\ldots,m$) satisfying homogeneous Dirichlet or Neumann boundary conditions on $\partial\Og\times(0,T)$. The system consists of $m$ equations 
of the form
\beqno{ieqn}(u_i)_t=\Div(\ag_i(\hat{u}_i)Du_i+\sum_{j< i}\bg_{ij}(\hat{u}_i)Du_j)+\eg_i|Du_i|^2+\sum_{j<i}c_{ij}|Du_j|^2+\mg_i,\eeq
where we denote $\hat{u}_i=(u_1,\ldots,u_i)$ and assume that $\eg_i\sup_Q|u_i|$ is small for all $i=1,\ldots,m$.

\bcoro{trisysthm} Consider the system \mref{fullsys0a} of $m$ equations of the form \mref{ieqn}. If $\mg_i\in L^\infty_{loc}(Q)$ then any bounded weak solution has bounded derivatives.
\ecoro

\bibliographystyle{plain}

\end{document}